\documentstyle[12pt]{article}
\def\Bbb R{{\rm \bf R}}
\def\proclaim#1{\vskip2mm{\bf #1}\em}
\def\endproclaim{\em \vskip2mm}
\def\tag#1{\eqno(#1)}
\def\gathered{\begin{array}{c}}
\def\endgathered{\end{array}}
\def\text{\mbox}

\begin{document}

\title {
The probe and enclosure methods for inverse obstacle scattering problems.
The past and present.}
\author{Masaru IKEHATA
\\
Department of Mathematics\\
Graduate School of Engineering\\
Gunma University, Kiryu 376-8515, JAPAN}
\date{ }
\maketitle



\tableofcontents

\section{The probe method for inverse obstacle scattering problems at a fixed wave number}
In this paper we consider inverse problems for partial
differential equations. We restrict ourself to the reconstruction
issue of the problems and refer the reader to \cite{ISA2} for
several aspects and uniqueness results in inverse problems for
partial differential equations.

More than ten years ago Ikehata discovered two methods for the
purpose of extracting information about the location and shape of
unknown {\it discontinuity} embedded in a known background medium
from observation data.  The methods are called the {\it probe} and
{\it enclosure} methods. This paper presents their past and recent
applications to inverse obstacle scattering problems of acoustic
wave.

The probe method was originally introduced in 1997 and published in \cite{IP}. 
Since then the method has been
applied to several inverse problems for partial differential
equations \cite{IPR, IP2, IPCOM, IPC, IPT} and still now some new
knowledge on the method itself added in \cite{IP3, IP4}.

In this section we present one of typical applications of the probe method published in \cite{IP2}.
Therein the author considered an inverse obstacle scattering problem at {\it a fixed wave number}.
We denote by $D$ and $B_R$ an unknown obstacle in $\Bbb R^3$
and open ball radius $R$, respectively.
We assume that: $D$ is an open set with smooth boundary satisfying $\overline D\subset B_R$
and that $B_R\setminus\overline D$ is connected.
$\partial B_R$ indicates the location of the emitters and the receivers.

Let $k>0$.  Given $y\in\partial B_R$ let $\Phi(x)=\Phi_D(x,y;k), x\in\Bbb R^3\setminus\overline D$
denote the solution of the problem:
$$\displaystyle
(\triangle +k^2)\Phi+\delta(\,\cdot\,-y)=0\,\,\text{in}\,\Bbb R^3\setminus\overline D,\,\,
\frac{\partial\Phi}{\partial\nu}=0\,\,\text{on}\,\partial D
$$
and the outgoing Sommerfeld radiation condition
$\displaystyle\lim_{r\longrightarrow\infty}r(\partial\Phi/\partial\nu-ik\Phi)=0$,
where $r=\vert x\vert$ and $\nu$ is the outward normal relative to $D$.

\noindent
{\bf Inverse Problem 1.1.}  {\it Fix} $k$.  Reconstruct $D$ from the surface data $\Phi_D(x,y;k)$ given at all
$x\in\partial B_R$ and $y\in\partial B_R$.

The $\Phi_D$ has the form $\displaystyle\Phi_D(x,y;k)=\Phi_0(x,y;k)+E_D(x,y;k)$,
where $E(x)=E_D(x,y;k)$ satisfies
$$\displaystyle
(\triangle+k^2)E=0\,\,\text{in}\,\Bbb R^3\setminus\overline D,\,\,
\frac{\partial E}{\partial\nu}=-\frac{\partial\Phi_0}{\partial\nu}\,\,\text{on}\,\partial D
$$
and the outgoing Sommerfeld radiation condition
$\displaystyle\lim_{r\longrightarrow\infty}r(\partial E/\partial\nu-ikE)=0$;
$\displaystyle\Phi_0(x,y;k)=e^{ik\vert x-y\vert}/(4\pi\vert x-y\vert)$.
The $E_D(x,y;k)$ is called the {\it scattered} wave field generated by the
{\it point source} $\delta(\,\cdot\,-y)$ located at $y$.  $\Phi_D(x,y;k)$ is called the {\it total} wave field.

In \cite{IP2} the author has established the following result.

\proclaim{\noindent Theorem 1.1.}
Assume that $k^2$ is not a Dirichlet eigenvalue for $-\triangle$ on $B_R$
nor an eigenvalue for $-\triangle$ on $B_R\setminus\overline D$ with
homogeneous Dirichlet boundary condition on $\partial B_R$ and Neumann boundary condition on $\partial D$.
Then one can reconstruct $D$ from $\Phi_D(x,y;k)$ given at {\it all} $x\in\partial B_R$ and $y\in\partial B_R$.
\endproclaim

A brief outline of the proof is as follows.
Set $\Omega=B_R$.
We starts with introducing two Dirichlet-to-Neumann maps for the Helmholtz equation in
$\Omega\setminus\overline D$ and $\Omega$.

Given $f\in H^{1/2}(\partial\Omega)$ let $u\in H^1(\Omega\setminus\overline D)$
be the weak solution of the elliptic problem
$$\begin{array}{c}
\displaystyle
(\triangle+k^2)u=0\,\,\text{in}\,\Omega\setminus\overline D,\,\,
\frac{\partial u}{\partial\nu}=0\,\,\text{on}\, \partial D,\,\,
u=f\,\,\text{on}\,\partial\Omega.
\end{array}
\tag {1.1}
$$
The map
$\displaystyle
\displaystyle\Lambda_D:f\longmapsto\partial u/\partial\nu\vert_{\partial\Omega}$
is called the {\it Dirichlet-to-Neumann map} associated with the elliptic problem.
Set also $\Lambda_D=\Lambda_0$ for $D=\emptyset$.

Theorem 1.1 is divided into two steps.

\proclaim{\noindent Step 1.}   One can calculate $\Lambda_0-\Lambda_D$ from
$\Phi_D(x,y;k)$ given at all $x\in\partial\Omega$ and $y\in\partial\Omega$.
\endproclaim

\proclaim{\noindent Step 2.}  One can reconstruct $D$ itself from the integral
$\displaystyle
\int_{\partial\Omega}(\Lambda_0-\Lambda_D)f\cdot\overline{f}dS$
for {\it infinitely many} $f$s independent of $D$.
\endproclaim
Note that the integral in Step 2 has the form
$$\displaystyle
\int_{\partial\Omega}(\Lambda_0-\Lambda_D)f\cdot\overline{f}dS
=\int_{\partial\Omega}\left(\frac{\partial v}{\partial\nu} \overline u-\frac{\partial u}{\partial\nu}\overline v\right)dS
$$
where $v=v(x)$, $x\in\Omega$ solves
$\displaystyle
(\triangle+k^2)v=0\,\,\text{in}\,\Omega,\,\,
v=f\,\,\text{on}\,\partial\Omega$; $u=u(x)$, $x\in\Omega\setminus\overline D$ solves (1.1) with $f=v\vert_{\partial\Omega}$.
Thus {\it infinitely many} $f$ means {\it infinitely many} $v$.

The step 1 consists of two parts.

\noindent
(i) Given $f$ find the solutions $g$ and $h$ of the integral equations
$$
\displaystyle
\int_{\partial\Omega}\Phi_0(x,y;k)g(y)dS(y)=f(x),\,\,
\int_{\partial\Omega}\Phi_D(x,y;k)h(y)dS(y)=f(x),\,\,x\in\partial\Omega.
$$

\noindent
(ii) Compute $(\Lambda_0-\Lambda_D)f$ by using solutions $g$ and $h$ in (i) by the formula
$\displaystyle (\Lambda_0-\Lambda_D)f=g-h$.

\noindent
This type of procedure, like (i) and (ii) have been known for the stationary Schr\"odinger equation \cite{N1}
and the proof  is an adaptation of the argument.
Thus the point is Step 2.

\subsection{Step 2.}

In this subsection we explain Step 2. Instead of the original
formulation of the probe method we employ a new one developed in
\cite{IP3, IP4}.

\subsubsection{Needle, Needle sequence}

{\bf\noindent Definition 1.1.}
Given a point $x\in\Omega$ we say that a non self-intersecting piecewise linear
curve $\sigma$ in $\overline\Omega$ is a {\it needle with tip at $x$}
if $\sigma$ connects a point on $\partial\Omega$
with $x$ and other points of $\sigma$ are contained in $\Omega$.
We denote by $N_x$ the set of all needles with tip at $x$.

Let $\mbox{\boldmath $b$}$ be a nonzero vector in $\Bbb R^3$.
Given $x\in\Bbb R^3$, $\rho>0$ and $\theta\in]0,\pi[$ set
$\displaystyle
C_x(\mbox{\boldmath $b$},\theta/2)
=\{y\in\,\Bbb R^3\,\vert\,
(y-x)\cdot\mbox{\boldmath $b$}>\vert y-x\vert\vert\mbox{\boldmath
$b$}\vert \cos(\theta/2)\}$
and $\displaystyle B_{\rho}(x)=\{y\in\,\Bbb R^3\,\vert\,\vert y-x\vert<\rho\}$.
A set having the form $\displaystyle
V=B_{\rho}(x)\cap C_x(\mbox{\boldmath $b$},\theta/2)$
for some $\rho$, $\mbox{\boldmath $b$}$, $\theta$ and $x$
is called a {\it finite cone} with {\it vertex} at $x$.

Let $G(y)$ be a solution of the Helmholtz equation in
$\Bbb R^3\setminus\{0\}$ such that, for any finite cone $V$ with vertex
at $0$
$$\displaystyle
\int_{V}\vert\nabla G(y)\vert^2 dy=\infty.
$$
Hereafter we {\it fix} this $G$.

{\bf\noindent Definition 1.2.}
Let $\sigma\in N_x$.  We call the sequence $\{v_n\}$
of $H^1(\Omega)$ solutions of the Helmholtz equation
a {\it needle sequence} for $(x,\sigma)$ if it satisfies,
for any compact set $K$ of $\Bbb R^3$ with $K\subset\Omega\setminus
\sigma$
$$\displaystyle
\lim_{n\longrightarrow\infty}(\Vert v_n(\,\cdot\,)-G(\,\cdot\,-x)\Vert_{L^2(K)}
+\Vert\nabla\{v_n(\,\cdot\,)-G(\,\cdot\,-x)\}\Vert_{L^2(K)})=0.
$$

The {\it existence} of the needle sequence is a consequence of
the {\it Runge approximation property} (cf.\cite{L}) for the Helmholtz equation
under the assumption on $k$:
$k^2$ is not a Dirichlet eigenvalue for $-\triangle$ on $\Omega$.
See the appendix of \cite{IP2} and A.1.Remark in the appendix of \cite{IP3} for the proof.
The {\it unique continuation property} of the solution of the Helmholtz equation
is essential.

\subsubsection{Special behaviour of the needle sequence}
In the following we do not assume that $k^2$ is not an eigenvalue for $-\triangle$
in $\Omega$ with Dirichlet boundary condition.

\proclaim{\noindent Lemma 1.1.}
Let $x\in\Omega$ be an arbitrary point and $\sigma\in N_x$.
Let $\{v_n\}$ be an arbitrary needle sequence for $(x,\sigma)$.
Then, for any finite cone $V$ with vertex at $x$ we have
$\displaystyle\Vert\nabla v_n\Vert_{L^2(V\cap\Omega)}\longrightarrow\infty$
as $n\longrightarrow\infty$.
\endproclaim

\proclaim{\noindent Lemma 1.2.}
Let $x\in\Omega$ be an arbitrary point and $\sigma\in N_x$.
Let $\{v_n\}$ be an arbitrary needle sequence for $(x,\sigma)$.
Then for any point $z\in\sigma$ and open ball $B$ centered at $z$ we have
$\displaystyle\Vert\nabla v_n\Vert_{L^2(B\cap\Omega)}\longrightarrow
\infty$ as $n\longrightarrow\infty$.
\endproclaim

\noindent
Note that from Definition 1.2 and Lemmas 1.1 and 1.2 one can recover $\sigma\in N_x$ itself
from the behaviour of any needle sequence for $(x,\sigma)$.

Summing up, we see that $\{v_n\}$ has two different sides:

\noindent
(A)  {\it converges} to singular solution $G(y-x)$ with singularity at $y=x$ outside $\sigma$;

\noindent
(B) {\it blows up} on $\sigma$.

\noindent
These different sides of needle sequences yield two sides of the probe method which we call Side A and Side B.

\subsubsection{Indicator function and Side A of the probe method}

Let $v$ satisfy $(\triangle+k^2)v=0$ in $\Omega$ and $u$ solve (1.1) with $f=v\vert_{\partial\Omega}$.
Set $\displaystyle w=u-v\,\,\text{in}\, \Omega\setminus\overline D$.
The $w$ satisfies
$$\displaystyle
(\triangle+k^2)w=0\,\,\text{in}\,\Omega\setminus\overline D,\,\,
w=0\,\,\text{on}\,\partial\Omega,\,\,
\frac{\partial w}{\partial\nu}=-\frac{\partial v}{\partial\nu}\,\,\text{on}\,\partial D.
\tag {1.2}
$$
Integration by parts yields
$$\begin{array}{c}
\displaystyle
\int_{\partial\Omega}(\Lambda_{\emptyset}-\Lambda_D)
(v\vert_{\partial\Omega})
\cdot\overline{v}dS
=\int_D\vert\nabla v\vert^2 dy-k^2\int_D\vert v\vert^2 dy\\
\\
\displaystyle
+\int_{\Omega\setminus\overline D}\vert\nabla w\vert^2dy
-k^2\int_{\Omega\setminus\overline D}\vert w\vert^2dy.
\end{array}
\tag {1.3}
$$

This motivates

{\bf\noindent Definition 1.3.}
The {\it indicator function} $I(x)$, $x\in\Omega\setminus\overline D$ is defined by the formula
$$\begin{array}{c}
\displaystyle
I(x)=\int_D\vert\nabla G(y-x)\vert^2 dy-k^2\int_D\vert G(y-x)\vert^2dy
+\int_{\Omega\setminus\overline D}
\vert\nabla w_x\vert^2 dy
-k^2\int_{\Omega\setminus\overline D}\vert w_x\vert^2dy,
\end{array}
$$
where $w_x$ is the unique weak solution of the problem:
$$\displaystyle
(\triangle+k^2)w=0\,\,\text{in}\,\Omega\setminus\overline D,\,\,
\frac{\partial w}{\partial\nu}=-\frac{\partial}{\partial\nu}(G(\,\cdot\,-x))\,\,\text{on}\,D,\,\,
w=0\,\,\text{on}\,\partial\Omega.
$$
The function $w_x$ is called the {\it reflected solution} by $D$.

The following theorem is based on the convergence property of needle sequences
and says that

\noindent
$\bullet$ one can calculate the value of the indicator function at an arbitrary point
outside $D$ from $\Lambda_0-\Lambda_D$;

\noindent
$\bullet$ the indicator function can not be continued across $\partial D$
as a bounded function in the whole domain.

Thus one can reconstruct $\partial D$ as the {\it singularity} of the field $I(x)$ which can be computed from the data
with needles and needle sequences.  That is the meaning of the following result.

\proclaim{\noindent Theorem A.}
It holds that

\noindent
$\bullet$  (A.1)  given $x\in\Omega\setminus\overline D$ and needle $\sigma$ with tip at $x$
if $\sigma\cap\overline D=\emptyset$,
then for any needle sequence $\{v_n\}$ for $(x,\sigma)$ we have
$\displaystyle
I(x)=\lim_{n\longrightarrow\infty}
\int_{\partial\Omega}(\Lambda_0-\Lambda_D)(v_n\vert_{\partial\Omega})\cdot\overline{v_n}dS$;

\noindent
$\bullet$  (A.2)  for each $\epsilon>0$
$\displaystyle
\sup\,\{I(x)\,\vert\,\displaystyle\text{dist}\,(x,\,D)>\epsilon\}<\infty$;

\noindent
$\bullet$  (A.3)  for any point $a\in\,\partial D$
$\displaystyle\lim_{x\longrightarrow a}I(x)=\infty$.

\endproclaim
The key for (A.3) is to establish
$\displaystyle
\limsup_{x\longrightarrow a}\Vert w_x\Vert_{L^2(\Omega\setminus\overline D)}<\infty$.
An outline of the proof is as follows.
Using the solution of the boundary value problem:
$\displaystyle
(\triangle +k^2)p=w_x$
in $\Omega\setminus\overline D$, $p=0$ on $\partial\Omega$
and $\partial p/\partial\nu=0$ on $\partial D$,
we have the expression
$$\displaystyle
\int_{\Omega\setminus\overline D}\vert w_x\vert^2dy
=\int_{\partial D}(p(x)-p(y))\overline{\frac{\partial\Phi_0}{\partial\nu}(y-x)}dS(y)
+k^2p(x)\int_D\overline{\Phi_0(y-x)}dy.
$$
Applying a standard regularity estimate of $p$:
$\Vert p\Vert_{H^2(\Omega\setminus\overline D)}\le C\Vert w_x\Vert_{L^2(\Omega\setminus\overline D)}$
and the Sobolev imbedding: $\vert p(x)-p(y)\vert\le C\vert x-y\vert^{1/2}\Vert p\Vert_{H^2(\Omega\setminus\overline D)}$,
$x,y\in\Omega\setminus\overline D$
and $\Vert p\Vert_{L^{\infty}(\Omega\setminus\overline D)}
\le C\Vert p\Vert_{H^2(\Omega\setminus\overline D)}$ to this right-hand side,
one gets an upper bound of $\Vert w_x\Vert_{L^2(\Omega\setminus\overline D)}$
which involves integrals of weakly singular kernels over $\partial D$ and $D$.

\subsection{Remark I.  Side B of the probe method and an open problem}

Since mathematically Theorem A is enough for establishing a
reconstruction formula, in the previous applications of the probe
method we did not consider the following natural question.

\noindent
$\bullet$ Let $x\in\Omega$ and $\sigma\in N_x$.  Let $\xi=\{v_n\}$ be a needle sequence for $(x,\sigma)$.
What happens on the sequence
$$\displaystyle
I(x,\sigma,\xi)_n\equiv\int_{\partial\Omega}(\Lambda_0-\Lambda_D)(v_n\vert_{\partial\Omega})\cdot\overline{v_n} dS, n=1,2,\cdots
$$
when $x$ is just located on the boundary of obstacles, inside or passing
through the obstacles?  We call sequence $\{I(x,\sigma,\xi)_n\}$
the {\it indicator sequence} for $(x,\sigma)$ and $\xi$.

In practice the tip of the needle can not move forward with
infinitely small step
and therefore in the scanning process with needle there is a
possibility of skipping the unknown boundary of obstacles, entering
inside or passing through obstacles. So for the
practical use of the probe method we have to clarify the behaviour
of the indicator sequence in those cases.  The answer to this
question is

\proclaim{\noindent Theorem B.}  Assume that $k^2$ is sufficiently small
(not specify here).
Let $x\in\Omega$ and $\sigma\in N_x$.
If $x\in\Omega\setminus\overline D$ and
$\sigma\cap D\not=\emptyset$
or $x\in\overline D$, then for any needle sequence $\xi=\{v_n\}$ for $(x,\sigma)$
we have $\displaystyle
\lim_{n\longrightarrow\infty}I(x,\sigma,\xi)_n=\infty$.
\endproclaim

\noindent
In the proof the blowing up property of needle sequences is essential.

{\it A sketch of the proof.}  For simplicity, we consider here
only a {\it single} obstacle case. We make use of two well known
Poincar\'e's inequalities:

(I) $\Vert w\Vert_{L^2(\Omega\setminus\overline D)}^2
\le C(\Omega\setminus\overline D)\Vert\nabla w\Vert_{L^2(\Omega\setminus\overline D)}^2$
for all $w\in H^1(\Omega\setminus\overline D)$ with $w=0$ on $\partial\Omega$;

(II) $\Vert v-v_D\Vert_{L^2(D)}^2\le C(D)\Vert\nabla v\Vert_{L^2(D)}^2$ for all $v\in H^1(D)$, where
$\displaystyle v_D=\int_D v dy/\vert D\vert$.

Let $A$ be an arbitrary Lebesgue measurable set with $A\subset D$, $\vert A\vert>0$ and $v\in L^2(D)$.
A simple argument in \cite{SO} gives
$\displaystyle
\Vert v-v_A\Vert_{L^2(D)}^2\le 2K_A\Vert v-v_D\Vert_{L^2(D)}^2$,
where $\displaystyle v_A=\int_A vdy/\vert A\vert$ and $K_A=1+\vert D\vert/\vert A\vert$.
A combination of this and (II) yields
$$\begin{array}{c}
\displaystyle
\int_{D}\vert v\vert^2 dy
\le
4K_A C(D)\int_{D}\vert\nabla v\vert^2 dy
+2\vert D\vert\vert v_A\vert^2.
\end{array}
\tag {1.4}
$$

Let $u=u_n$ solve (1.1) with $f=v_n\vert_{\partial\Omega}$ and set $w_n=u_n-v_n$.
It follows from (1.3), (I) and (1.4) that
$$\begin{array}{c}
\displaystyle
I(x,\sigma,\xi)_n
\ge (1-k^2C(\Omega\setminus\overline D))\int_{\Omega\setminus\overline D}\vert\nabla w_n\vert^2 dy\\
\\
\displaystyle
+(1-4k^2K_AC(D))\int_{D}\vert\nabla v_n\vert^2 dy
-2k^2\vert D\vert\vert (v_n)_A\vert^2.
\end{array}
$$
Thus if $k$ satisfies $\displaystyle k^2C(\Omega\setminus\overline D)\le 1$, then we have
$$
\displaystyle
I(x,\sigma,\xi)_n
\ge
(1-4k^2K_AC(D))\int_{D}\vert\nabla v_n\vert^2 dy
-2k^2\vert D\vert\vert (v_n)_A\vert^2.
$$
Write $\displaystyle
1-4k^2K_AC(D)=1-8k^2C(D)-4k^2(K_A-2)C(D)$.
Here we make $k$ smaller in such a way that $8k^2C(D)<1$.
Using an exhaustion of $\Omega\setminus\sigma$, one can construct
$A\subset D$ in such a way that $
\vert A\vert\approx\vert D\vert$ and $\overline A\subset\Omega\setminus\sigma$.
Since $K_A-2=\vert D\vert/\vert A\vert-1$, one gets $1-4k^2K_AC(D)>0$.
Note also that the sequence $\{(v_n)_A\}$ is always {\it convergent} for a fixed $A$.
Thus the blowing up property of the indicator sequence is governed by that of
the sequence $\{\Vert\nabla v_n\Vert_{L^2(D)}^2\}$.

A combination of Theorems A and B yields another characterization of
the obstacle.

\proclaim{\noindent Corollary 1.1.}
Assume the smallness of $k^2$ same as Theorem B.
A point $x\in\,\Omega$ belongs to $\Omega\setminus\overline D$ if and only if
there exists a needle $\sigma$ with tip at $x$ and needle sequence $\xi$
for $(x,\sigma)$ such that the indicator sequence
is bounded from above.
\endproclaim

\noindent
Needless to say, this automatically gives a uniqueness theorem, too.

An open problem in the foundation of the probe method is the following.

\noindent
{\bf Open problem 1.1.}
Can one remove the smallness of $k^2$ in Theorem B?

Here are some closely related technical questions.

\noindent
$\bullet$ Is it true ?: if $x\in\Omega\setminus\overline D$ and
$\sigma\cap D\not=\emptyset$
or $x\in\overline D$, then
$$\displaystyle
\lim_{n\longrightarrow\infty}\frac
{\displaystyle
\Vert v_n\Vert_{L^2(D)}}
{\displaystyle\Vert\nabla v_n\Vert_{L^2(D)}}
=0.
\tag {1.5}
$$

\noindent
$\bullet$  Let $u=u_n$ solve (1.1) with $f=v_n\vert_{\partial\Omega}$ and set $w_n=u_n-v_n$.
We know that if $x\in\overline D$, then $\Vert\nabla w_n\Vert_{L^2(\Omega\setminus\overline D)}\longrightarrow\infty$
as $n\longrightarrow\infty$ (\cite{IPT}).  
The question is: identify the points in $\overline\Omega\setminus D$ that really contribute
the blowing up of $\nabla w_n$.  See \cite{IP3} for an example in the case when $k=0$.

\noindent
$\bullet$  Is it true ?: if $x\in\Omega\setminus\overline D$ and
$\sigma\cap D\not=\emptyset$
or $x\in\overline D$, then
$$\displaystyle
\lim_{n\longrightarrow\infty}\frac{\displaystyle
\Vert w_n\Vert_{L^2(\Omega\setminus\overline D)}}
{\displaystyle
\Vert\nabla v_n\Vert_{L^2(D)}}=0.
$$
See \cite{IP3, IPC, IPT} for more information on these questions.

\subsection{Remark II.  An explicit needle sequence}

From Lemmas 1.1 and 1.2 we know that given $\sigma\in N_x$  the energy of an arbitrary needle sequence  $\{v_n\}$ for $(x,\sigma)$
blows up on $\sigma$.
However, it will be difficult to understand the behaviour of $v_n(y)$ at each $y\in\sigma$.
In this subsection, we give a family of special solutions of the Helmholtz equation with two parameters
that yields an explicit needle sequence for a straight needle.  We call such a family a generator of needle sequence.

The contents of this subsection are based on the classical materials developed by
Yarmukhamedov, Mittag-Leffler and Vekua.

\subsubsection{Yarmukhamedov}

The following fact is taken from the article \cite{Y}.

\proclaim{\noindent Theorem 1.2.}
Let $K(w)$ be an entire function such that: $K(w)$ is real for real $w$;
$K(0)=1$; for each $R>0$  and $m=0, 1, 2$
$\displaystyle
\sup_{\vert\text{Re}\,w\vert<R}\vert K^{(m)}(w)\vert
<\infty$.

\noindent
Define
$$\displaystyle
-2\pi^2\Phi_K(x)
=\int_0^{\infty}\text{Im}\,\left(\frac{K(w)}{w}\right)\frac{du}
{\sqrt{\vert x'\vert^2+u^2}},
$$
where $w=x_3+i\sqrt{\vert x'\vert^2+u^2}$ and $x'=(x_1,x_2)\not=(0,0)$.
Then one has the expression
$\displaystyle\Phi_K(x)=1/(4\pi\vert x\vert)+H_K(x)$
where $H_K$ satisfies $\displaystyle
\triangle H_K(x)=0\,\,\text{in}\,\Bbb R^3$.

\endproclaim
Note that $\Phi_K$ can be identified with a unique distribution in the whole space
and satisfies $\displaystyle\triangle\Phi_K(x)+\delta(x)=0\,\,\text{in}\,\Bbb R^3$.

\noindent
{\bf Example 1.} $K(w)\equiv 1$.
In this case we have $\displaystyle\Phi_K(x)=1/(4\pi\vert x\vert)$.
This is because of
$$\displaystyle
\frac{1}{4\pi\vert x\vert}
=\int_{-\infty}^{\infty}\frac{du}{4\pi^2(\vert x\vert^2+u^2)}
$$
and
$$\displaystyle
\frac{1}{\vert x\vert^2+u^2}
=-\text{Im}\,\left(\frac{1}{x_3+i\sqrt{\vert x'\vert^2+u^2}}\right)
\frac{1}{\sqrt{\vert x'\vert^2+u^2}}.
$$
Thus for general $K$ we have
$$\displaystyle
H_K(x)
=-\frac{1}{2\pi^2}
\int_0^{\infty}\text{Im}\,\left(\frac{K(w)-1}{w}\right)\frac{du}
{\sqrt{\vert x'\vert^2+u^2}}.
$$

\noindent
{\bf Example 2.} $K(w)\equiv e^{\tau w}$. $\tau>0$ a parameter.
In \cite{ISLAB} the author pointed out that
$\displaystyle\Phi_K(x)$ with this $K$
coincides with the {\it Faddeev Green function} $G_z(x)$
with $z=\tau(e_3+ie_1)$:
$$\displaystyle
G_z(x)=\frac{e^{x\cdot z}}{(2\pi)^3}
\int_{\Bbb R^3}\frac{e^{ix\cdot\eta}}{\vert\eta\vert^2-i2z\cdot\eta}d\eta.
$$
The Faddeev Green function has been applied to
several inverse boundary value/scattering problems by Sylvester-Uhlmann \cite{SU},
Novikov \cite{No}, Nachman \cite{N1}, et al..

\subsubsection{Mittag-Leffler}

Let $0<\alpha\le 1$.
The entire function of the complex variable $w$
$$\displaystyle
E_{\alpha}(w)=1+\frac{w}{\Gamma(1+\alpha)}
+\frac{w^2}{\Gamma(1+2\alpha)}+\frac{w^3}{\Gamma(1+3\alpha)}+\cdots,\,
$$
is introduced in \cite{Mi} and called the {\it Mittag-Leffler function}.

It is known that $K(w)=E_{\alpha}(\tau w)$ with $\tau>0$ satisfies the condition
in Theorem 1.2 (cf. \cite{B}).
In \cite{YE} Yarmukhamedov applied this function with a fixed $\alpha$ to the Cauchy problem
for the Laplace equation in two dimensions.

\subsubsection{Vekua}

The {\it Vekua
transform} $v\longmapsto T_{k}v$ in three dimensions \cite{V} takes
the form
$$\displaystyle
T_{k}v(y)=v(y)-\frac{k\vert y\vert}{2}\,\int_0^1v(ty)J_1(k\vert y\vert\sqrt{1-t})
\,\sqrt{\frac{t}{1-t}}\,dt
$$
where $J_1$ stands for the Bessel function of order $1$ of the first kind.

The important property of this transform
is: if $v$ is harmonic in the whole space, then $T_{k}v$ is a solution of the Helmholtz
equation $\triangle u+k^2 u=0$ in the whole space.

\subsubsection{Generator of needle sequence}

Using materials introduced by Yarmukhamedov, Mittag-Leffler and Vekua, the author
found an explicit needle sequence when the needle is given by a {\it segment}.

Given $0<\alpha\le 1$ and $\tau>0$ define $v(y;\alpha,\tau)=-H_K(y)$, $y\in\Bbb R^3$,
where $K(w)\equiv E_{\alpha}(\tau w)$.  This $v$ is harmonic in the whole space
and thus the function $v^{k}(y;\alpha,\tau)=T_{k}v(y;\alpha,\tau)$, $y\in\Bbb R^3$
satisfies the Helmholtz equation in the whole space.

\proclaim{\noindent Theorem 1.3(\cite{IP4}).}  Let $x\in\Omega$ and $\sigma$
be a straight needle with tip at $x$ directed to
$\omega=(0,0,1)^T$, that means: $\sigma$ has the expression $\sigma=
\{x+s\,\omega\,\vert\,0\le s\le l\}$ with $l>0$. 
Then the function
$v^{k}(\,\cdot\,-x;\alpha,\tau)\vert_{\Omega}$ as $\alpha\longrightarrow 0$ and
$\tau\longrightarrow\infty$
generates a needle sequence for $(x,\sigma)$ with $G=G_{k}$ given by
$$\displaystyle
G_{k}(y)=\text{Re}\,\left(\frac{\displaystyle
e^{ik\vert y\vert}}{
\displaystyle
4\pi\vert y\vert}\right).
$$
\endproclaim
Note that since the function
$$\displaystyle
\frac{\displaystyle\sin\,k\vert y\vert}{\displaystyle
4\pi\vert y\vert},\,y\in\Bbb R^3
$$
satisfies the Helmholtz equation in the whole space, the function
$$\displaystyle
v^k(y-x;\alpha,\tau)+i\,\frac{\displaystyle\sin\,k\vert y-x\vert}{\displaystyle
4\pi\vert y-x\vert},\,\,y\in\Omega
$$
generates also a needle sequence for $(x,\sigma)$ with
$$\displaystyle
G(y)=\frac{\displaystyle
e^{ik\vert y\vert}}{\displaystyle
4\pi\vert y\vert}.
\tag {1.6}
$$

Thus now we have an
explicit generator of a needle sequence for a straight needle with (1.6).  This makes the probe method
completely {\it explicit} in the case when one uses only such a needle.
Everything is reduced to the choice of small $\alpha$ and large $\tau$.

This is very important also in the {\it singular sources
method} by Potthast \cite{P1} since in his method one has to construct
the density of the Herglotz wave function (cf.
\cite{CK}) that approximates locally the fundamental solution of
the Helmholtz equation in a domain like $\Omega\setminus\sigma$.
However, Theorem 1.3 shows that instead
one can consider only a simpler problem: construct the density of
the Herglotz wave function that approximates
$v^k(y-x;\alpha,\tau)$ on the whole {\it boundary} of a geometrically simpler domain like
a {\it ball}.

\noindent
{\bf Open problem 1.2.}
It would be interesting: do the numerical testing of the probe and singular sources methods
in three dimensions with this {\it explicit}
needle sequence.

\noindent
{\bf Open problem 1.3.}
A mathematically interesting question is: find a generator of a needle sequence
for a general needle.

Note that Yarmukhamedov \cite{YM} made use of $\Phi_K(y-x)$ itself not its regular part $H_K(y-x)$
to give a Carleman function which yields a representation of the solution
of the Cauchy problem for the Laplace equation in three dimensions.

Finally we give a remark that is closely related to Open problem 1.1.
In \cite{IP4} an explicit formula of the precise values of $v^{k}(y-x;\alpha,\tau)$
on the line $y=x+\,s\,\omega (-\infty<s<\infty)$ is given.
They are:

\noindent
$\bullet$ if $y=x+s\,\omega$ with $s\not=0$, then
$$\begin{array}{c}
\displaystyle
v^{k}(y-x;\alpha,\tau)=
\frac{1}{4\pi}
\frac{\displaystyle
E_{\alpha}(\tau s)-\cos\,ks}
{s}
-\frac{k}{4\pi}
\int_0^1(1-w^2)^{-1/2}E_{\alpha}(\tau(1-w^2)s)
J_1(ksw)dw;
\end{array}
$$

\noindent
$\bullet$ if $y=x$, then
$$\displaystyle
v^{k}(y-x;\alpha,\tau)\vert_{y=x}
=\frac{\tau}{4\pi\Gamma(1+\alpha)}.
$$

\noindent
Moreover, we see that
$\nabla v^{k}(y-x;\alpha,\tau)$ on the line
$y=x+s\,\omega\,(-\infty<s<\infty)$ is parallel to $\omega$.  In particular,
we have
$$\displaystyle
\nabla v^{k}(y-x;\alpha,\tau)\vert_{y=x}
=\frac{\tau^2}{4\pi\,\Gamma(1+2\alpha)}\,\omega.
$$
It seems that the behaviuor of $v^{k}(y-x;\alpha,\tau)$ and its gradient at $y=x$ suggest
the validity of (1.5).

\section{The enclosure method for inverse obstacle scattering problems at a fixed wave number}

The {\it enclosure method} was introduced by the author in
\cite{IE} and has been applied to several inverse problems for
partial differential equations. In this section we present its
applications to inverse obstacle scattering problems at a fixed
wave number.

\subsection{The enclosure method with infinitely many data}

The method applied to inverse obstacle scattering problems
is based on the asymptotic behaviour of the function (we call the indicator function again)
$$\displaystyle
\tau\longmapsto
\int_{\partial\Omega}(\Lambda_0-\Lambda_D)(v\vert_{\partial\Omega})\cdot\overline{v}dS,
$$
where $\displaystyle
v=e^{x\cdot(\tau\,\omega+i\sqrt{\tau^2+k^2}\,\omega^{\perp})}$
having large parameter $\tau$; both $\omega$ and
$\omega^{\perp}$ are unit vectors and perpendicular to each
other.

This $v$ satisfies the Helmholtz equation $\triangle
v+k^2v=0$ in the whole space and divides the whole space into two
parts: if $x\cdot\omega>t$, then $e^{-\tau t}\vert
v\vert\longrightarrow \infty$ as $\tau\longrightarrow\infty$;
if $x\cdot\omega<t$, then $e^{-\tau t}\vert v\vert\longrightarrow
0$ as $\tau\longrightarrow\infty$.

The method yielded the convex hull of unknown {\it sound-soft} obstacles
by checking the behaviour of the indicator function.
It virtually checks whether given $t$ the half space $x\cdot\omega>t$ touches unknown obstacles.

In \cite{IE2} an extraction formula of an {\it sound-hard
obstacle} $D\subset\Bbb R^3$ with a constrained on the {\it
Gaussian curvature} of $\partial D$ from Dirichlet-to-Neumann map
$\Lambda_D$ has been established. Its precise statement rewritten
with the present style is the following.

Let us recall the {\it support function} of $D$:
$\displaystyle
h_D(\omega)=\sup_{x\in\,D}x\cdot\omega,\,\,\omega\in S^2$.
The convex hull of $D$ is given by the set
$\displaystyle\cap_{\omega\in S^2}\{x\in\Bbb  R^3\,\vert\,x\cdot\omega<h_D(\omega)\}$.
Therefore, knowing $h_D(\omega)$ for a $\omega$ yields an estimation of the convex hull of $D$
from above.
\proclaim{\noindent Theorem 2.1.}
Assume that the set $\displaystyle
\{x\in\partial D\,\vert\,x\cdot\omega=h_D(\omega)\}$
consists of only one point and the Gaussian curvature of $\partial D$
doesn't vanish at the point.
Then the formula
$$\displaystyle
\lim_{\tau\longrightarrow\infty}\frac{1}{2\tau}
\log\left\vert\int_{\partial\Omega}(\Lambda_{0}-\Lambda_D)
(v\vert_{\partial\Omega})
\cdot\overline{v}dS\right\vert
=h_D(\omega),
$$
is valid.  Moreover, we have:

if $t>h_D(\omega)$, then
$$\displaystyle
\lim_{\tau\longrightarrow\infty}
\int_{\partial\Omega}(\Lambda_{0}-\Lambda_D)
(e^{-\tau t}v\vert_{\partial\Omega})
\cdot\overline{e^{-\tau t}v}dS=0;
$$

if $t<h_D(\omega)$, then
$$\displaystyle
\lim_{\tau\longrightarrow\infty}
\int_{\partial\Omega}(\Lambda_{0}-\Lambda_D)
(e^{-\tau t}v\vert_{\partial\Omega})
\cdot\overline{e^{-\tau t}v}dS=\infty;
$$

if $t=h_D(\omega)$, then
$$\displaystyle
\liminf_{\tau\longrightarrow\infty}
\int_{\partial\Omega}(\Lambda_{0}-\Lambda_D)
(e^{-\tau t}v\vert_{\partial\Omega})
\cdot\overline{e^{-\tau t}v}dS>0.
$$
\endproclaim

Note that:  if one considers the Dirichlet boundary condition
$\text{$u=0$ on $\partial D$}$
instead of the Neumann boundary condition $\partial u/\partial\nu=0$
on $\partial D$, one can {\it drop} the assumption on $\omega$ and the Gaussian
curvature of $\partial D$.  See \cite{IE} for this result.
Thus we propose

\noindent
{\bf Open problem 2.1.} Remove the curvature condition in Theorem 2.1.

{\it A sketch of the proof of Theorem 2.1.} Let $u$ solve (1.1) with $f=v\vert_{\partial\Omega}$
and set $\displaystyle
w=u-v\,\,\text{in}\, \Omega\setminus\overline D$.
The $w$ satisfies (1.2).
We have three lemmas. \proclaim{\noindent Lemma 2.1.} There exists a
positive constat $C(k)$ such that for all $\omega\in S^2$,
$\tau>0$
$$\begin{array}{c}
\displaystyle
2\tau^2\int_De^{2\tau x\cdot\omega}dx- k^2\int_{\Omega\setminus\overline D}
\vert w\vert^2dx
\le
\int_{\partial\Omega}(\Lambda_{0}-\Lambda_D)
(v\vert_{\partial\Omega})
\cdot\overline{v}dS
\le C(k)(\tau^2+k^2)\int_De^{2\tau x\cdot\omega}dx.
\end{array}
$$
\endproclaim

This is a consequence of the representation formula (1.3)
and the estimate
$\displaystyle
\Vert w\Vert_{H^1(\Omega\setminus\overline D)}
\le C(k)\Vert v\Vert_{H^1(D)}$.

\proclaim{\noindent Lemma 2.2.}
$$\displaystyle
\liminf_{\tau\longrightarrow\infty}e^{-2\tau h_D(\omega)}\tau^2
\int_De^{2\tau x\cdot\omega}dx>0.
$$

\endproclaim
The proof of this lemma can be done by slicing $D$ with the planes
$x\cdot\omega=h_D(\omega)-s$ with $0<s<<1$.

\proclaim{\noindent Lemma 2.3.}
Assume that the set
$\displaystyle
\{x\in\partial D\,\vert\,x\cdot\omega=h_D(\omega)\}$
consists of the only one point and the Gaussian curvature of $\partial D$
doesn't vanish at the point.
Then
$$\displaystyle
\lim_{\tau\longrightarrow\infty}
\frac{\displaystyle\int_{\Omega\setminus\overline D}\vert w\vert^2dx}
{\displaystyle 2\tau^2\int_De^{2\tau x\cdot\omega}dx}
=0.
$$
\endproclaim
From Lemmas 2.1, 2.2 and 2.3 one knows that there exist positive constants $C_1$, $C_2$ and $\tau_0>0$ such that
for all $\tau\ge\tau_0$
$$\displaystyle
C_1e^{2\tau h_D(\omega)}\le
\int_{\partial\Omega}(\Lambda_{0}-\Lambda_D)
(v\vert_{\partial\Omega})
\cdot\overline{v}dS
\le C_2\tau^2e^{2\tau h_D(\omega)}.
$$
All the statements in Theorem 2.1 now follows from these
estimates.

Finally we describe the outline of the proof of Lemma 2.3.
One can find $p\in H^2(\Omega\setminus\overline D)$
such that
$\displaystyle
(\triangle+k^2)p=\overline w$ in $\Omega\setminus\overline D$,
$p=0$ on $\partial\Omega$ and
$\partial p/\partial\nu=0$ on $\partial D$.
From the Sobolev imbedding and the estimate
$\displaystyle
\Vert p\Vert_{H^2(\Omega\setminus\overline D)}\le C(k)
\Vert w\Vert_{L^2(\Omega\setminus\overline D)}$
we have:
$\displaystyle
\vert p(x)-p(y)\vert\le C(k)\vert x-y\vert^{1/2}\Vert w\Vert_{L^2(\Omega\setminus\overline D)}$
and
$\displaystyle
\sup_{x\in\overline\Omega\setminus D}\vert p(x)\vert\le C(k)
\Vert w\Vert_{L^2(\Omega\setminus\overline D)}$.

Let $x_0$ be the point in the set
$\{x\in\partial D\,\vert\,x\cdot\omega=h_D(\omega)\}$.
Since
$\displaystyle
\int_{\partial D}(\partial v/\partial\nu)dS(x)=
-k^2\int_D v dx$,
one can write
$$\begin{array}{c}
\displaystyle
\int_{\Omega\setminus\overline D}\vert w\vert^2 dx
=-\int_{\partial D}p\frac{\partial v}{\partial\nu}dS(x)
=\int_{\partial D}\{p(x_0)-p(x)\}\frac{\partial v}{\partial\nu}dS(x)
+k^2p(x_0)\int_D v dx.
\end{array}
$$
From these one gets
$$\begin{array}{c}
\displaystyle
\int_{\Omega\setminus\overline D}\vert w\vert^2 dx
\le C(k)\left(\sqrt{2\tau^2+k^2}\int_{\partial D}\vert x_0-x\vert^{1/2}
e^{\tau x\cdot\omega}dS(x)
+\int_De^{\tau x\cdot\omega}dx\right)\Vert w\Vert_{L^2(\Omega\setminus\overline D)}
\end{array}
$$
and this thus yields
$$\begin{array}{c}
\displaystyle
\int_{\Omega\setminus\overline D}\vert w\vert^2 dx
\le C(k)
\left\{\left(\tau\int_{\partial D}\vert x_0-x\vert^{1/2}
e^{\tau x\cdot\omega}dS(x)\right)^2
+\left(\int_De^{\tau x\cdot\omega}dx\right)^2\right\}.
\end{array}
$$
The Schwarz inequality yields
$$\displaystyle
\left(\int_De^{\tau x\cdot\omega}dx\right)^2
\le\vert D\vert\int_De^{2\tau x\cdot\omega}dx.
$$
Thus from this and Lemma 2.2 one knows that
it suffices to prove
$$\displaystyle
\lim_{\tau\longrightarrow\infty}\tau e^{-\tau h_D(\omega)}\int_{\partial D}\vert x_0-x\vert^{1/2}
e^{\tau x\cdot\omega}dS(x)=0.
$$
In fact, one gets
$$\displaystyle
\tau e^{-\tau h_D(\omega)}\int_{\partial D}\vert x_0-x\vert^{1/2}
e^{\tau x\cdot\omega}dS(x)
=O(\tau^{-1/4}).
$$
This is proved by using a localization at $x_0$ and a local coordinates at the point.

\subsection{The enclosure method with a single incident plane wave}

The idea started with considering an inverse boundary value problem for the Laplace equation
in two dimensions in \cite{IE3}.
Five years later in \cite{IE4} the idea was applied to an inverse obstacle scattering problem in two dimensions.
The problem is to reconstruct a two dimensional obstacle from the Cauchy data on a circle surrounding the obstacle of the total wave field 
for a {\it single} incident plane wave with a {\it fixed} wave number.

In this subsection we assume that
 $D$ is {\it polygonal}, that is,
$D$ takes the form $D_1\cup\cdots\cup D_m$ with $1\le m<\infty$
where each $D_j$ is open and a polygon; $\overline D_j\cap\overline D_{j'}=\emptyset$ if $j\not=j'$.

The total wave field $u$ outside obstacle $D$ satisfies
$$
\displaystyle
\triangle u+k^2 u=0\,\,\text{in}\,\Bbb R^2\setminus\overline D,\,\,
\frac{\partial u}{\partial\nu}=0\,\,\text{on}\,\partial D
$$
and the scattered wave $w=u-e^{ikx\cdot d}$  with $k>0$ and $d\in S^1$ satisfies
the outgoing Sommerefeld radiation condition
$\displaystyle
\lim_{r\longrightarrow\infty}
\sqrt{r}(\partial w/\partial r-ikw)=0$,
where $r=\vert x\vert$.

Let $B_R$ be an open disc with radius $R$ satisfying $\overline
D\subset B_R$.  We assume that $B_R$ is {\it known}.  Our data are  $u$ and $\partial u/\partial\nu$ on $\partial B_R$.
Let $\omega$ and $\omega^{\perp}$ be two unit
vectors perpendicular to each other.  Set
$z=\tau\omega+i\sqrt{\tau^2+k^2}\omega^{\perp}$ with $\tau>0$ and
$v(x;z)=e^{x\cdot z}$.  Recall $h_D(\omega)=\sup_{x\in D}x\cdot\omega$.

\proclaim{\noindent Theorem 2.2.}
Assume that the set
$\displaystyle
\partial D\cap\{x\in\Bbb R^2\,\vert\,x\cdot\omega=h_D(\omega)\}$
consists of only one point.
Then the formula
$$\displaystyle
\lim_{\tau\longrightarrow\infty}\frac{1}{\tau}
\log\left\vert
\int_{\partial B_R}\left(\frac{\partial u}{\partial\nu}v(x;z)
-\frac{\partial v}{\partial\nu}(x;z)u\right)
dS(x)
\right\vert
=h_D(\omega),
$$
is valid.  Moreover, we have:

if $t\ge h_D(\omega)$, then
$$\displaystyle
\lim_{\tau\longrightarrow\infty}
\left\vert
\int_{\partial B_R}\left(\frac{\partial u}{\partial\nu}e^{-\tau t}v(x;z)
-e^{-\tau t}\frac{\partial v}{\partial\nu}(x;z)u\right)
dS(x)
\right\vert=0;
$$

if $t<h_D(\omega)$, then
$$\displaystyle
\lim_{\tau\longrightarrow\infty}\left\vert
\int_{\partial B_R}\left(\frac{\partial u}{\partial\nu}e^{-\tau t}v(x;z)
-e^{-\tau t}\frac{\partial v}{\partial\nu}(x;z)u\right)
dS(x)
\right\vert
=\infty.
$$
\endproclaim

{\it Sketch of the proof.} The one of key points is: {\it
introducing a new parameter} $s$ instead of $\tau$ by the equation
$\displaystyle s=\sqrt{\tau^2+k^2}+\tau$, we obtain, as
$s\longrightarrow\infty$ the complete asymptotic expansion
$$\displaystyle
\int_{\partial B_R}\left(\frac{\partial u}{\partial\nu}v(x;z)
-\frac{\partial v}{\partial\nu}(x;z)u\right)
dS(x)e^{-i\sqrt{\tau^2+k^2}x_0\cdot\omega^{\perp}-\tau h_D(\omega)}
\sim
-i\sum_{n=2}^{\infty}\frac{e^{i\frac{\pi}{2}\lambda_n}k^{\lambda_n}\alpha_nK_n}
{s^{\lambda_n}}.
\tag {2.1}
$$
Here the $\lambda_n$ describes the {\it singularity} of $u$ at a {\it corner} and in this case explicitly given by
the formula $\displaystyle\lambda_n=(n-1)\pi/\Theta$,
where $\Theta$ denotes the outside angle of $D$ at $x_0\in\partial
D\cap\{x\in\Bbb R^2\,\vert\,x\cdot\omega=h_D(\omega)\}$ and thus
satisfies $\pi<\Theta<2\pi$; $K_n$ are constants depending on
$\lambda_n$, $\omega$ and shape of $D$ around $x_0$; $\alpha_2$,
$\alpha_3$, $\cdots$ are the coefficients of the convergent series expansion of $u$
with polar coordinates at a corner:
$$\displaystyle
u(r,\theta)=\alpha_1J_0(kr)+\sum_{n=2}^{\infty}\alpha_n J_{\lambda_n}(kr)\cos \lambda_n\theta,\,\,0<r<<1, 0<\theta<\Theta.
$$

Now all the statements in Theorem 2.2 follow from (2.1) and
another key point: $\exists n\ge 2$\,\, $\alpha_nK_n\not=0$.
This is due to a contradiction argument.
Assume that the assertion is not true, that is,
$\forall n\ge 2$\,\,$\alpha_nK_n=0$.

First we consider the case when $\Theta/\pi$ is {\it irrational}.
In this case we see that $\forall n\ge 2$\,\,$K_n\not=0$.  Thus
$\alpha_n=0$ and this yields $u(r,\theta)=\alpha_1J_0(kr)$ near a
corner.  Since this right-hand side is an entire solution of the
Helmholtz equation, the unique continuation property of the
solution of the Helmholtz equation yields $u(x)=\alpha_1J_0(k\vert
x-x_0\vert)$ in $\Bbb R^2\setminus\overline D$. However, we see
that the asymptotic behaviour of this right-hand and left-hand
sides are completely different.  Contradiction.

Next consider the case when $\Theta/\pi$ is a {\it rational}.  By
carefully checking the constant $K_n$ we know that for each $n\ge
2$ with $K_n=0$ the $\lambda_n$ becomes an {\it integer}. From the
assumption of the contradiction argument one knows if $n$
satisfies $K_n\not=0$, then $C_n=0$.  Thus we have the expansion
$$\displaystyle
u(r,\theta)=\sum_{n_j}C_{n_j}J_{\lambda_{n_j}}(kr)\cos\,\lambda_{n_j}\theta,
$$
where $n_j\ge 2$ satisfy $K_{n_j}=0$.  Since $\lambda_{n_j}$ is an
integer and $\lambda_{n_j}\Theta=(n_j-1)\pi$, from this right-hand
side one gets: for all $r$ with $0<r<<1$ $\displaystyle\partial
u/\partial\theta(r,\pi)=\partial
u/\partial\theta(r,\Theta-\pi)=0$. Then a reflection argument
(\cite{AD}) yields that this is true for all $r>0$.  However, from
this together with the asymptotic behaviour of $\nabla u$ one can
conclude that incident direction $d$ has to be parallel to two
linearly independent vectors which are directed along the lines
$\theta=\pi$ and $\theta=\Theta-\pi$.  Contradiction.

Remarks are in order.

\noindent
$\bullet$  In Theorem 2.2 one uses the Cauchy data on the circle surrounding the
obstacle as the observation data.  However, $\partial u/\partial\nu$ on $B_R$ can be
calculated from $u$ on $\partial B_R$ by solving an exterior Dirichlet problem for
the Helmholz equation.

\noindent
$\bullet$  In \cite{IHER} a similar formula has been established by using
the {\it far field
pattern} $F_D(\varphi,d;k)$, $\varphi\in S^1$ of scattered wave $w=u-e^{ikx\cdot d}$ for fixed $d$ and $k$
which determines the leading term of the asymptotic expansion of $w$ at infinity in the following sense
$$\displaystyle
w(r\varphi)\sim\frac{\displaystyle
e^{ikr}}{\displaystyle
\sqrt{r}}F_D(\varphi,d;k)\,\,r\longrightarrow\infty.
$$
Moreover, therein instead of volumetric obstacle, similar formulae for {\it thin} sound-hard
obstacle (or {\it screen} )
also have been established with {\it two} incident plane waves.

\noindent
$\bullet$  In \cite{OI} the numerical testing of a method based on results in \cite{ILAY2,
IE4, IHER} has been reported.

\noindent
$\bullet$  It would be interesting to consider the case when the total wave
$u$ satisfies the equation $\nabla\cdot\gamma\nabla u+k^2 u=0$ in
$\Bbb R^2$ where $\gamma(x)=1$ for $x\in\Bbb R^2\setminus D$ and
$\gamma(x)=A_j$ for $x\in D_j$, $j=1,\cdots, m$; each $A_j$ are
positive constants and $A_j\not=1$. The author thinks that this
case becomes extremely difficult because of the complicated
behaviour of $u$ at a corner.  However
we propose

\noindent
{\bf Open problem 2.2.}  Establish Theorem 2.2 for $u$ above.

\noindent
See \cite{IE100} for $k=0$ and \cite{ITRANS} for the equation
$\nabla\cdot\gamma\nabla u+k^2\,\gamma u=0$.

\noindent
$\bullet$  For recent applications of the enclosure method with a single
measurement for a system arising in linear theory of elasticity we
have \cite{IH1, IH2, IH}.  However, their extension to the elastic
wave with a single incident plane wave remains open.  It is a
challenging problem to be solved.

\section{Inverse obstacle scattering problems with dynamical data over a finite
time interval}

Previously we considered only the stationary or time harmonic problem. In this
section we consider how one can use the data over a {\it finite
time interval} to extract information about the location and shape
of unknown obstacles. In \cite{ISA, R, R2} some uniqueness results
have been established, however, it seems that mathematically
rigorous study of the reconstruction issue in this type of problem
has not been paid much attention. Note that: there are some results
\cite{LP, M, MT} in the context of the Lax-Phillips scattering
theory, which give the convex hull of an unknown obstacle,
however, the data are taken from $t=0$ to $t=\infty$.

The purpose of this section is to introduce a new and simple
method in \cite{IW} which is an application of the idea developed
in \cite{IHEAT,IK} and employs the data over a finite time
interval on a known surface surrounding unknown obstacles.

\subsection{New development of the enclosure method}

In order to explain the basic idea, in this subsection we present
an application to the one-space dimensional wave equation which is
taken from Appendix B in \cite{IHEAT}. Let $a>0$ and $c>0$. Let
$u=u(x,t)$ be a solution of the problem:
$$\begin{array}{c}
\displaystyle
\frac{1}{c^2}u_{tt}=u_{xx}\,\,\text{in}\,]0,\,a[\times]0,\,T[,\,\,
cu_x(a,t)=0\,\,\text{for}\,t\in\,]0,\,T[,\\
\\
\displaystyle
u(x,0)=0,\,\,
u_t(x,0)=0\,\,\text{in}\,]0,\,a[.
\end{array}
$$
The quantity $c$ denotes the propagation speed
of the signal governed by the equation.

\noindent
{\bf Inverse Problem 3.1.}
Assume that $a$ is {\it unknown}.
Extract $a$ from $u(0,t)$ and $u_x(0,t)$
for $0<t<T$.

\proclaim{\noindent Theorem 3.1.}
Let $u_x(0,t)\in L^2(0,\,T)$ satisfy the
condition: there exists a real number $\mu$ such that
$$\displaystyle
\liminf_{\tau\longrightarrow\infty}\tau^{\mu}\left\vert\int_0^Tu_x(0,t)e^{-\tau t}dt\right\vert>0.
\tag {3.1}
$$
Let $T>2a/c$ and $\displaystyle
v(x,t)=v(x,t;\tau)=e^{-\tau(x/c+t)}$.
Then the formula
$$\displaystyle
\lim_{\tau\longrightarrow\infty}
\frac{1}{\tau}
\log\left\vert
\int_0^T\left(-cv_x(0,t)u(0,t)+cu_x(0,t)v(0,t)\right)dt
\right\vert
=-2a/c,
\tag {3.2}
$$
is valid.

\endproclaim

Some remarks are in order.

\noindent
$\bullet$  The $v$ satisfies the wave equation $(1/c)^2v_{tt}=v_{xx}$ and satisfies: if $x+ct>0$, then
$v(x,t)\longrightarrow 0$ as $\tau\longrightarrow\infty$; if $x+ct<0$, then $v(x,t)\longrightarrow +\infty$
as $\tau\longrightarrow\infty$.

\noindent
$\bullet$  The quantity $2a/c$ coincides with the travel time of a
signal governed by the wave equation  with propagation speed $c$ which starts at the boundary
$x=0$ and initial time $t=0$, reflects another boundary $x=a$ and
returns to $x=0$.  Thus the restriction $T>2a/c$ is quite reasonable and
does not against the well known fact: the wave equation has
the {\it finite propagation property}.

\noindent
$\bullet$  The condition (3.1)
ensures that $u_x(0,t)$ can not be identically zero in an interval
$]0,T'[\subset\,]0,\,T[$. Therefore surely a signal occurs at the
initial time.  However, it should be emphasized that the formula
(3.2) makes use of the {\it averaged value} of the measured data
with an {\it exponential weight} over the observation time.
This is a completely different idea from the well known approach in
nondestructive evaluation by sound wave: monitoring of the first
{\it arrival time} of the {\it echo}, one knows the travel time.

{\it A sketch of the proof of Theorem 3.1.}
Introduce the function $w$ by the formula
$$\displaystyle
w(x)=w(x;\tau)=\int_0^Tu(x,t)e^{-\tau t}dt,\,\,0<x<a.
$$
It holds that
$$\begin{array}{c}
\displaystyle
c^2w^{''}-\tau^2w=e^{-\tau T}(u_t(x,T)+\tau u(x,T))\,\,\text{in}\,]0,\,a[,\,\,
cw'(a)=0.
\end{array}
$$
Then, this together with integration by parts gives the expression
$$\begin{array}{c}
\displaystyle
e^{2a\tau/c}\int_0^T\left(-cv_x(0,t)u(0,t)+cu_x(0,t)v(0,t)\right)dt\\
\\
\displaystyle
=\tau w(a)e^{a\tau/c}
-c^{-1}e^{-\tau(T-(2a/c))}\int_0^a(u_t(\xi,T)+\tau u(\xi,T))e^{-\xi\tau/c}d\xi.
\end{array}
$$
Now (3.2) can be checked by studying the asymptotic behaviour of this right-hand side
with the help of the expression
$$\begin{array}{l}
\displaystyle
w(a)
=-\frac{2cw'(0)}
{\displaystyle
\tau\left(e^{a\tau/c}-e^{-a\tau/c}\right)}-\frac{e^{-\tau T}}
{\displaystyle
\tau\left(e^{a\tau/c}-e^{-a\tau/c}\right)}
\\
\\
\displaystyle
\times
\left\{\int_0^a(u_t(\xi,T)+\tau u(\xi,T))e^{-\xi\tau/c}d\xi
+\int_0^a(u_t(\xi,T)+\tau u(\xi,T))e^{\xi\tau/c}d\xi\right\}
\end{array}
$$
together with (3.1).

The proof presented here heavily relies on the spaciality of
one-space dimension. In \cite{IK} we found another method for the
proof which works also for higher space dimensions and applied it
to a similar problem for the heat equation.  In the following two
subsections we present further applications of the method to the
wave equations.

\subsection{Sound-hard obstacle}
Let $D\subset\Bbb R^3$ be a bounded open set with smooth boundary
such that $\Bbb R^3\setminus\overline D$ is connected.
Denote by $\nu$ the unit outward normal to $\partial D$.
Let $0<T<\infty$.

Given $f\in L^2(\Bbb R^3)$ with compact support satisfying
$\text{supp}\,f\cap\overline D=\emptyset$ let $u=u(x,t)$
satisfy the initial boundary value problem:
$$\begin{array}{c}
\displaystyle
\partial_t^2u-\triangle u=0\,\,\text{in}\,(\Bbb R^3\setminus\overline D)\times\,]0,\,T[,\,\,
\frac{\partial u}{\partial\nu}=0\,\,\text{on}\,\partial D\times\,]0,\,T[,\\
\\
\displaystyle
u(x,0)=0,\,\,
\partial_tu(x,0)=f(x)\,\,\text{in}\,\Bbb R^3\setminus\overline D.
\end{array}
$$
Let $\Omega$ be a bounded domain with smooth boundary such that $\overline D\subset\Omega$
and $\Bbb R^3\setminus\overline\Omega$ is connected.
Denote by the same symbol $\nu$ the unit outward normal to $\partial\Omega$.

The $\partial\Omega$ is considered as the location of the receivers
of the acoustic wave produced by an emitter located at the support of $f$.
In this section we consider the following problem.

{\bf\noindent Inverse Problem 3.2.}
Assume that $D$ is {\it unknown.}
Extract information about the location and shape of $D$
from $u$ on $\partial\Omega\times]0,\,\,T[$ for some fixed {\it known} $f$
satisfying $\text{supp}\,f\cap\overline\Omega=\emptyset$ and $T<\infty$.

\noindent
Note that $u$ in $(\Bbb R^3\setminus\overline\Omega)\times]0,\,\,T[$ can be computed from
$u$ on $\partial\Omega\times\,]0,\,T[$ by
the formula
$$\displaystyle
u=z\,\,\text{in}\,(\Bbb R^3\setminus\overline\Omega)\times\,]0,\,T[
\tag {3.3}
$$
where $z$ solves the initial boundary value problem in $\Bbb R^3\setminus\overline\Omega$:
$$\begin{array}{c}
\displaystyle
\partial_t^2z-\triangle z=0\,\,\text{in}\,(\Bbb R^3\setminus\overline\Omega)\times\,]0,\,T[,\,\,
z=u\,\,\text{on}\,\partial\Omega\times\,]0,\,T[,\\
\\
\displaystyle
z(x,0)=0,\,\,
\partial_tz(x,0)=f(x)\,\,\text{in}\,\Bbb R^3\setminus\overline\Omega.
\end{array}
$$

Thus the problem can be reformulated as

{\bf\noindent Inverse Problem 3.2'.}
Extract information about the location and shape of $D$ from
$u$ in $(\Bbb R^3\setminus\overline\Omega)\times]0,\,T[$ for some known $f$ satisfying
$\text{supp}\,f\cap\overline\Omega=\emptyset$
and $T<\infty$.

Now we state the result.
Let $B$ be an open ball with $\overline B\cap\overline\Omega=\emptyset$.
Choose the initial data $f\in L^2(\Bbb R^3)$ in such a way that:

\noindent
(I1) $f(x)=0$ a.e. $x\in\Bbb R^3\setminus B$;

\noindent
(I2) there exists a positive constant $C$ such that $f(x)\ge C$ a.e. $x\in B$
or $-f(x)\ge C$ a.e. $x\in B$.

Set
$$\displaystyle
w(x;\tau)=\int_0^T e^{-\tau t}u(x,t)dt,\,\,x\in\Bbb R^3\setminus\overline D,\,\,\tau>0.
$$
Our result is the following extraction formula from
$w$ and $\partial w/\partial\nu$ on $\partial\Omega\times\,]0\,\,T[$ which
can be computed from the data $u$ in $(\Bbb R^3\setminus\overline\Omega)\times\,]0,\,T[$.
\proclaim{\noindent Theorem 3.2.}
Let $\tau>0$ and $v\in H^1(\Bbb R^3)$ be the weak solution of
$$\displaystyle
(\triangle-\tau^2)v+f(x)=0\,\,\text{in}\,\Bbb R^3.
\tag {3.4}
$$
If the observation time $T$ satisfies
$$\displaystyle
T>2\text{dist}\,(D,B)-\text{dist}\,(\Omega,B),
\tag {3.5}
$$
then there exists a $\tau_0>0$ such that, for all $\tau\ge\tau_0$
$$\displaystyle
\int_{\partial\Omega}\left(\frac{\partial v}{\partial\nu}w-\frac{\partial w}{\partial\nu}v\right)dS>0
$$
and the formula
$$\displaystyle
\lim_{\tau\longrightarrow\infty}
\frac{1}{2\tau}\log
\int_{\partial\Omega}\left(\frac{\partial v}{\partial\nu}w-\frac{\partial w}{\partial\nu}v\right)dS
=-\text{dist}\,(D,B),
\tag {3.6}
$$
is valid.
\endproclaim

Some remarks are in order.

\noindent
$\bullet$  The $v$ is unique and is given by the explicit form
$$\displaystyle
v(x;\tau)=\frac{1}{4\pi}\int_B
\frac{e^{-\tau \vert x-y\vert}}{\vert x-y\vert}f(y)dy,\,\,x\in\Bbb R^3.
$$

\noindent
$\bullet$  The quantity $\text{dist}\,(D,B)+\sqrt{\vert\partial B\vert/4\pi}$ coincides with the distance from the center of $B$ to $D$
and thus (3.6) yields the information about $d_D(p)$ for a given point $p$ in $\Bbb R^3\setminus\overline\Omega$.

\noindent
$\bullet$  It is easy to see that $2\text{dist}\,(D,B)-\text{dist}\,(\Omega,B)\ge
l(\partial B,\partial D,\partial\Omega)$, where
$\displaystyle
l(\partial B,\partial D,\partial\Omega)
=\inf\,\{\vert x-y\vert+\vert y-z\vert\,\vert\,x\in\partial B\,, y\in\partial D,\, z\in\partial\Omega\}$.
This is the {\it minimum length of the broken paths} that start at $x\in\partial B$ and reflect at $y\in\partial D$ and return
to $z\in\partial\Omega$.
Therefore (3.5) ensures that $T$ is greater than the {\it first arrival time}
of a signal with the unit propagation speed that starts at a point on $\partial B$ at $t=0$,
reflects at a point on $\partial D$ and goes to a point on $\partial\Omega$.

The main part of the proof of Theorem 3.2 is to show that
$$\displaystyle
\liminf_{\tau\longrightarrow\infty}\tau^{4}e^{2\tau\,\text{dist}\,(D,B)}
\int_{\partial\Omega}\left(\frac{\partial v}{\partial\nu}w-\frac{\partial w}{\partial\nu}v\right)dS>0.
\tag {3.7}
$$
It is a consequence of the following representation formula which corresponds to (1.3) and the estimate for $v$:
$$\begin{array}{c}
\displaystyle
\int_{\partial\Omega}\left(\frac{\partial v}{\partial\nu}w-\frac{\partial w}{\partial\nu}v\right)dS\\
\\
\displaystyle
=\int_D\vert\nabla v\vert^2dx+\tau^2\int_D\vert v\vert^2dx
+\int_{\Bbb R^3\setminus\overline D}\vert\nabla(w-v)\vert^2dx
+\tau^2\int_{\Bbb R^3\setminus\overline D}\vert w-v\vert^2dx
\\
\\
\displaystyle
+e^{-\tau T}\int_{\Bbb R^3\setminus\overline D}(w-v)(\partial_tu(x,T)+\tau u(x,T))dx
-e^{-\tau T}\int_{\Omega\setminus\overline D}
(\partial_tu(x,T)+\tau u(x,T))vdx;
\end{array}
$$
$$\displaystyle
\liminf_{\tau\longrightarrow\infty}\tau^{6}e^{2\tau\text{dist}\,(D,B)}\int_D\vert v\vert^2 dx>0.
\tag {3.8}
$$
Note that the precise values of $4$ and $6$ of $\tau^{4}$ in (3.7) and $\tau^{6}$ in (3.8), respectively
are not essential.

\subsection{Penetrable obstacle}

The method in the former subsection can be applied to a more general case.
Given $f\in L^2(\Bbb R^3)$ with compact support
let $u=u(x,t)$
satisfy the initial value problem:
$$\begin{array}{c}
\displaystyle
\partial_t^2u-\nabla\cdot\gamma\nabla u=0\,\,\text{in}\,\Bbb R^3\times\,]0,\,T[,\\
\\
\displaystyle
u(x,0)=0,\,\,
\partial_tu(x,0)=f(x)\,\,\text{in}\,\Bbb R^3,
\end{array}
\tag {3.9}
$$
where $\gamma=\gamma(x)=(\gamma_{ij}(x))$ satisfies:  for each $i,j=1,2,3$ $\gamma_{ij}(x)=\gamma_{ji}(x)\in L^{\infty}(\Bbb R^3)$;
there exists a positive constant $C$ such that $\gamma(x)\xi\cdot\xi\ge C\vert\xi\vert^2$ for all $\xi\in\Bbb R^3$
and a. e. $x\in\Bbb R^3$.

We assume: there exists a bounded open set $D$ with a smooth boundary
such that $\gamma(x)$ a.e. $x\in\Bbb R^3\setminus D$ coincides with the $3\times 3$ identity matrix $I_3$.
Write $h(x)=\gamma(x)-I_3$ a.e. $x\in D$.

Our second inverse problem is the following.

{\bf\noindent Inverse Problem 3.3.}  Assume that both $D$ and $h$ are {\it unknown}
and that one of the following two conditions is satisfied:

\noindent
(A1)  there exists a positive constant $C$ such that $-h(x)\xi\cdot\xi\ge\vert\xi\vert^2$
for all $\xi\in\Bbb R^3$ and a.e. $x\in D$;

\noindent
(A2)  there exists a positive constant $C$ such that $h(x)\xi\cdot\xi\ge\vert\xi\vert^2$
for all $\xi\in\Bbb R^3$ and a.e. $x\in D$.

\noindent
Let $\Omega$ be a bounded domain with smooth boundary such that $\overline D\subset\Omega$.
Extract information about the location and shape of $D$
from $u$ on $\partial\Omega\times]0,\,\,T[$ for some fixed {\it known} $f$
satisfying $\text{supp}\,f\cap\overline\Omega=\emptyset$ and $T<\infty$.

Note that $u$ in $(\Bbb R^3\setminus\overline\Omega)\times]0,\,\,T[$ can be computed from
$u$ on $\partial\Omega\times]0,\,T[$ by the exactly same formula as (3.3)
and thus the problem can be reformulated again as

{\bf\noindent Inverse Problem 3.3'.}
Extract information about the location and shape of $D$ from
$u$ in $(\Bbb R^3\setminus\overline\Omega)\times\,]0,\,T[$ for some known $f$ satisfying
$\text{supp}\,f\cap\overline\Omega=\emptyset$ and $T<\infty$.

Now we state our second result.

\proclaim{\noindent Theorem 3.3.}
Assume that $\gamma$ satisfies (A1) or (A2).
Let $f$ satisfy (I1) and (I2) in subsection 3.2 and $v$ be the weak solution of (3.4).
Let $T$ satisfies (3.5)
and $w$ be given by
$$\displaystyle
w(x;\tau)=\int_0^T e^{-\tau t}u(x,t)dt,\,\,x\in\Bbb R^3,\,\,\tau>0
$$
with solution $u$ of (3.9).
If (A1) is satisfied, then
there exists a $\tau_0>0$ such that, for all $\tau\ge\tau_0$
$$\displaystyle
\int_{\partial\Omega}\left(\frac{\partial v}{\partial\nu}w-\frac{\partial w}{\partial\nu}v\right)dS>0;
$$
if (A2) is satisfied, then there exists a $\tau_0>0$ such that, for all $\tau\ge\tau_0$
$$\displaystyle
-\int_{\partial\Omega}\left(\frac{\partial v}{\partial\nu}w-\frac{\partial w}{\partial\nu}v\right)dS>0.
$$
In both cases we have
$$\displaystyle
\lim_{\tau\longrightarrow\infty}
\frac{1}{2\tau}\log
\left\vert\int_{\partial\Omega}\left(\frac{\partial v}{\partial\nu}w-\frac{\partial w}{\partial\nu}v\right)dS\right\vert
=-\text{dist}\,(D,B).
$$
\endproclaim
The key points for the proof are an estimate for $\nabla v$ similar to (3.8) and
the following two representation formula:
$$\begin{array}{c}
\displaystyle
\int_{\partial\Omega}
\{(\nabla v\cdot\nu)w-(\gamma\nabla w\cdot\nu)v\}dS
=-\int_D h\nabla v\cdot\nabla v dx\\
\\
\displaystyle
+\int_{\Bbb R^3}\gamma\nabla(w-v)\cdot\nabla(w-v)dx
+\tau^2\int_{\Bbb R^3}\vert w-v\vert^2dx\\
\\
\displaystyle
+e^{-\tau T}\int_{\Bbb R^3}(\partial_tu(x,T)+\tau u(x,T))(w-v)dx
-e^{-\tau T}
\int_{\Omega}(\partial_tu(x,T)+\tau u(x,T))vdx;
\end{array}
$$
$$\begin{array}{c}
\displaystyle
-\int_{\partial\Omega}\{(\nabla v\cdot\nu)w-(\gamma\nabla w\cdot\nu)v\}dS
=\int_Dh\nabla w\cdot\nabla wdx\\
\\
\displaystyle
+\int_{\Bbb R^3}\nabla(v-w)\cdot\nabla(v-w)dx
+\tau^2\int_{\Bbb R^3}\vert v-w\vert^2dx\\
\\
\displaystyle
-e^{-\tau T}\int_{\Bbb R^3}(\partial_tu(x,T)+\tau u(x,T))(v-w)dx
+e^{-\tau T}\int_{\Omega}(\partial_tu(x,T)+\tau u(x,T))vdx.
\end{array}
$$

\section{Summary and further research direction}

In this paper we presented: past applications of the probe and enclosure methods
to inverse obstacle scattering problems with a fixed wave number
and related open problems; recent applications of the enclosure method
to {\it inverse obstacle scattering problems with dynamical data over a finite time interval}.

In particular, in Section 3 we presented a new and simple method in \cite{IW} for a typical class
of inverse obstacle scattering problems
that employs the values of the wave field over a {\it finite} time interval on a known surface surrounding unknown obstacles
as the observation data.  The wave field is generated by an initial data localized outside the surface
and its form is not specified except for the condition on the support.
The method {\it explicitly} yields information
about the location and shape of the obstacles {\it more than} the convex hull.

It would be interesting to apply the method
presented in Section 3 to other time dependent problems in
electromagnetism(e.g., {\it subsurface radar} \cite{DGS}, {\it microwave tomography} \cite{S}),
linear elasticity, classical fluids etc..
Those applications belong to our future plan.

$$\quad$$

\centerline{{\bf Acknowledgements}}

This research was partially supported by Grant-in-Aid for
Scientific Research (C)(No. 21540162) of Japan  Society for the
Promotion of Science.

\vskip1cm
\noindent
e-mail address

ikehata@math.sci.gunma-u.ac.jp

\end{document}